\documentclass{amsart}
\usepackage{amssymb}
\newtheorem{thm}{Theorem}
\newtheorem{lem}{Lemma}

\begin{document}

\bibliographystyle{plain}

\title[Rado\v s Baki\'c]{On the representation of a polynomial as a sum of monomials}

\author[]{Milo\v s Arsenovi\'c}
\author[]{Rado\v s Baki\'c}

\address{Faculty of mathematics, University of Belgrade, Studentski Trg 16, 11000 Belgrade, Serbia}
\email{\rm arsenovic@matf.bg.ac.rs}

\address{U\v citeljski Fakultet, Beograd, Serbia}
\email{\rm bakicr@gmail.com}

\date{}

\begin{abstract}
We investigate systems of $m \geq n$ polynomial equations in $n$ complex variables, where polynomials involved are power sums, and give estimates of the solutions. These results are applied to obtain a representation of arbitrary polynomial of degree $n$ as the average of powers $(z-z_j)^n$.
\end{abstract}

\maketitle

\footnotetext[1]{Mathematics Subject Classification 2010 Primary 12E10 secondary 12E12.  Key words
and Phrases: Symmetric polynomials, estimates of solutions.}



Two basic families of symmetric polynomials in $n$ variables $x_1, \ldots, x_n$ are elementary symmetric polynomials
$\sigma_k = \sigma_k(x_1, \ldots, x_n)$ and power sums $s_k = s_k(x_1, \ldots, x_n)$ defined by
$$\sigma_k = \sum_{1 \leq i_1 < \cdots < i_k \leq n} x_{i_1}x_{i_2} \cdots x_{i_k}, \qquad 1 \leq k \leq n,$$

$$\sigma_k = 0, \qquad k > n,$$
and
$$ s_k = x_1^k + \cdots + x_n^k, \qquad k \geq 1.$$
These families are related by well-known identities:
\begin{equation}\label{epk}
k\sigma_k = \sum_{i=1}^k (-1)^{i-1} s_i \sigma_{k-i}, \qquad k \geq 1.
\end{equation}

These relations are linear with respect to both families $\sigma_j$ and $s_k$, hence we can express members of one family in
terms of members of the other family. Explicitly
\begin{equation}\label{ek} \sigma_k = \frac{1}{k!} \begin{vmatrix}
s_1 & 1 & 0 & \cdots & 0 \\
s_2 & s_1 & 2 & \cdots & 0 \\
: & : & : & \cdots & :\\
s_{k-1} & s_{k-2} & s_{k-3} & \cdots & k-1\\
s_k & s_{k-1} & s_{k-2} & \cdots & s_1
\end{vmatrix} = \Delta_k (s_1, \ldots, s_k)
\end{equation}
and
\begin{equation}\label{pk} s_k = \begin{vmatrix}
\sigma_1 & 1 & 0 & \cdots & 0\\
2\sigma_2 & \sigma_1 & 1 & \cdots & 0\\
: & : & : & \cdots & :\\
(k-1)\sigma_{k-1} & \sigma_{k-2} & \sigma_{k-3} & \cdots & 1\\
k\sigma_k & \sigma_{k-1} & \sigma_{k-2} & \cdots & \sigma_1
\end{vmatrix} = D_k (\sigma_1, \ldots, \sigma_k).
\end{equation}

Let us note that equations (\ref{ek}) and (\ref{pk}) are equivalent, both being equivalent to (\ref{epk}).

We say that a monomial $As_1^{a_1} \ldots s_k^{a_k}$ is $k$-regular in $s_1, \ldots, s_k$ if $a_1 + 2a_2 + \cdots +
ka_k = k$ and $A$ is a real number such that $\mbox{sgn} A = (-1)^{k+a_1 + \cdots + a_k}$.

\begin{lem}
For any $k \geq 1$ polynomial $\sigma_k$ can be expressed as a sum of $k$-regular monomials in $s_1, \ldots, s_k$.
\end{lem}

{\it Proof}. We use induction on $k$ and relations (\ref{epk}). For $k=1$  we have  $\sigma_1= s_1$ as desired. Now we assume the
statement is true for $1 \leq j < k$. Because of (\ref{epk}) and induction hypothesis we need only to prove that if  $As_1^{a_1} \ldots s_{k-i}^{a_{k-i}}$ is $(k-i)$-regular then  $(-1)^{i-1}s_iAs_1^{a_1} \ldots s_{k-i}^{a_{k-i}}$ is $k$-regular, which can be easily checked. $\Box$

Our estimates of the moduli of solutions of systems of polynomial equations are based on the
following theorem due to Tchakaloff, see \cite{Tc}.

\begin{thm}\label{thmtc}
Let $f_i(z)$   be a finite set of complex polynomials of degree $n$,  all having a positive leading coefficient. If all $f_i(z)$  have their zeros in the disk $|z-z_0| \leq r$ , then polynomial $F(z) = \sum_i f_i(z)$  has all its zeros in the disk $|z-z_0| \leq C_n r$, where
\begin{equation}\label{cn} C_n = \frac{1}{\sin \frac{\pi}{2n}}.
\end{equation}
\end{thm}

Further results related to the above theorem can be found in \cite{Ru}.

\begin{thm}\label{thmrb}
Suppose that the following system in variables $z_1, \ldots, z_n$  is given
\begin{equation}\label{system}
z_1^j + z_2^j + \cdots + z_n^j = b_j, \quad 1 \leq j \leq n.
\end{equation}
Then it has unique solutions for $z_1, \ldots, z_n$ (up to permutations). Also $|z_i| \leq C_nM$ where $C_n$ is the constant given
in (\ref{cn}) and $M = \displaystyle\max_{1 \leq j \leq n} |b_j|^{1/j}$.
\end{thm}

{\it Proof.} Our system can be written in the form $s_j(z_1, \ldots, z_n) = b_j$, $1 \leq j \leq n$. Using (\ref{ek}) and
(\ref{pk}), we can rewrite it in an equivalent form
\begin{equation}\label{form}
\sigma_k(z_1, \ldots, z_n) = c_k, \qquad 1 \leq k \leq n,
\end{equation}
where $c_k = \Delta_k(b_1, \ldots, b_k)$. By the Fundamental theorem of algebra and Viet's rule, the above system has unique solution  for $z_j$  (up to permutations). Since two previous systems are equivalent, the first part of our lemma is proved. Let $(z_1, \ldots, z_n)$
be a solution; it suffices, by symmetry, to prove $|z_n| \leq C_n M$. Let $s_k^\prime = z_1^k + \cdots + z_{n-1}^k$ be power sum polynomials in $z_1, \ldots, z_{n-1}$, and let $\sigma_k^\prime$ be the corresponding elementary symmetric polynomials.  By the above lemma, it follows that $\sigma_n^\prime$ can be expressed as a sum of $n$-regular monomials $P_\lambda$, $1 \leq \lambda \leq l$ in $s_1^\prime, \ldots, s_j^\prime$. Let $A_\lambda \prod_{j=1}^n (s_j^\prime)^{a_{j, \lambda}}$ be one of them.

Using (\ref{system}) we get $s_j^\prime = b_j - z_n^j$, hence we can rewrite previous monomial as
$A_\lambda \prod_{j=1}^n (b_j - z_n^j)^{a_{j, \lambda}}$. Hence, since our monomial is $n$-regular, it can be seen as a polynomial
$P_\lambda (z_n)$ in $z_n$  of degree $n$, the sign of the leading coefficient is equal to
$(-1)^{a_{1, \lambda} + \cdots  + a_{n, \lambda}} {\rm sgn} A_\lambda = (-1)^n$. However, $\sigma_n^\prime$  is identically equal to zero, hence we have equation $\sum_{\lambda = 1}^l P_\lambda (z_n) = 0$  where ${\rm deg} P_\lambda = n$ for all $1 \leq \lambda \leq l$ and all the leading coefficients of $P_\lambda$ have the same sign $(-1)^n$. Note that all such polynomials have all
their zeros in the disc $|z| \leq M$. Hence by Theorem \ref{thmtc} their sum has all its zeros in the disc $|z| \leq C_n M$, therefore
$|z_n| \leq MC_n$ as required. $\Box$

\begin{thm}\label{any}
Any complex polynomial
\begin{equation}\label{pol}
f(z) = z^n + a_{n-1}z^{n-1} + \ldots + a_1 z + a_0
\end{equation}
can be written in the form
\begin{equation}\label{ave}
f(z) = \frac{1}{n} \sum_{i=1}^n (z-z_i)^n,
\end{equation}
where $z_1, \ldots, z_n$  are uniquely determined up to their order. Also, if $M = \displaystyle\max_{1 \leq j \leq n} |c_j|^{1/j}$, where
$c_j = \frac{n a_{n-j}}{\binom{n}{j}}$ then we have an estimate $|z_i| \leq C_n M$, where $C_n$ is the constant given in (\ref{cn}).
\end{thm}

{\it Proof.} Relation (\ref{ave}) leads to the following system of equations:
\begin{equation}
z_1^j + \cdots + z_n^j = (-1)^j c_j, \qquad 1 \leq j \leq n.
\end{equation}
According   to Theorem \ref{thmtc} above system has desired solution which proves the theorem. $\Box$

Next, motivated by Theorem \ref{thmrb}, we consider the following system in $m$ complex variables, where $m \geq n$:
\begin{equation}\label{ttt}
\frac{1}{m} (z_1^j + \cdots + z_m^j) = A_j^j, \qquad 1 \leq j \leq n.
\end{equation}
Since $m \geq n$, it follows from Theorem \ref{thmrb} that a solution exists. We are interested in estimates of solutions.
A general problem is to find the best possible constant $K_{n,m}$ such that for any complex numbers $A_1, \ldots, A_n$ there is a solution to the system (\ref{ttt}) satisfying $\displaystyle\max_{1 \leq k \leq m} |z_k| \leq K_{n,m} A$, where $A =
\displaystyle\max_{1 \leq j \leq n} |A_j|$. A weaker version of the same problem is to find $K_n =
\displaystyle\inf_{m \geq n} K_{n,m}$. An asymptotic behavior of $K_n$ would be of interest as well.

Our last theorem gives a partial answer to the questions posed above. Let us define a sequence $D_n$ inductively by $D_1 = 1$ and
\begin{equation}\label{rec}
D_{n+1} = 1 + \left( 1 + D_n^{n+1} \right)^{1/n+1}, \qquad n \geq 1.
\end{equation}

\begin{thm}\label{nfact}
If $m = n!$ then the system (\ref{ttt}) has a solution $(z_k)_{1 \leq k \leq n!}$ such that
\begin{equation}
\max_{1 \leq k \leq n!} |z_k| \leq D_n \max_{1 \leq j \leq n} |A_j|.
\end{equation}
\end{thm}

{\it Proof.} We proceed by induction on $n$. The case $n=1$ is trivial, so we assume the statement is true for $n$. Let $A_1, \ldots,
A_n, A_{n+1}$ be arbitrary complex numbers and set $A = \displaystyle\max_{1 \leq j \leq n+1} |A_j|$. By inductive hypothesis there exist complex numbers $z_k$, $1 \leq k \leq n!$ such that
\begin{equation}\label{indh}
\frac{1}{n!} \sum_{k=1}^{n!} z_k^j = A_j^j, \qquad 1 \leq j \leq n
\end{equation}
and
\begin{equation}\label{estn}
\max_{1 \leq k \leq n!} |z_k| \leq D_n \max_{1 \leq j \leq n} |A_j| \leq D_n A.
\end{equation}

Set $\xi = \exp \frac{2\pi i}{n+1}$ and define $(n+1)!$ complex numbers
\begin{equation}\label{wkl}
w_{k,l} = z_k + \xi^l K, \qquad 1 \leq k \leq n!, \quad 0 \leq l \leq n
\end{equation}
where constant $K$ will be chosen later on. Using elementary identity
$$\sum_{l=0}^n \xi^{lr} = \begin{cases} n+1 & \text{if $n+1$ divides $r$}\\
0 & \text{otherwise} \end{cases}$$
we deduce that for every $1 \leq j \leq n$ and $1 \leq k \leq n!$ we have
\begin{align*}
\sum_{l=0}^n w_{k,l}^j & =  \sum_{l=0}^n (z_k + \xi^l K)^j  =  \sum_{l=0}^n \sum_{r=0}^j \binom{j}{r} z_k^{j-r} K^r \xi^{lr} \\
& =  \sum_{r=0}^n \binom{j}{r} z_k^{j-r} K^r \sum_{l=0}^n \xi^{lr} =  (n+1)z_k^j
\end{align*}

Therefore, using (\ref{indh}) we have
\begin{equation*}
\frac{1}{(n+1)!}  \sum_{k=1}^{n!} \sum_{l=0}^n w_{k,l}^j = \frac{1}{n!} \sum_{k=1}^n z_k^j = A_j^j
\end{equation*}
for every $1 \leq j \leq n$, which means that the $(n+1)!$ numbers $w_{k,l}$ satisfy the first $n$ equations of the system (\ref{ttt}),
with $m = (n+1)!$, {\it regardless} of the choice of the constant $K$. Next we have

\begin{align*}
\frac{1}{(n+1)!} \sum_{k,l}^{n+1} w_{k,l}^{n+1} & = \frac{1}{(n+1)!} \sum_{k=1}^{n!} \sum_{l=0}^n (z_k + \xi^l K)^{n+1} \\
& = \frac{1}{(n+1)!} \sum_{k=1}^{n!} \sum_{l=0}^n \sum_{r=0}^{n+1} \binom{n+1}{r} z_k^{n+1-r} K^r \xi^{lr} \\
& = \frac{1}{(n+1)!} \sum_{k=1}^{n!} \sum_{r=0}^{n+1} \binom{n+1}{r} z_k^{n+1-r} K^r \sum_{l=0}^n \xi^{lr} \\
& = \frac{1}{n!} \sum_{k=1}^{n!} (z_k^{n+1} + K^{n+1})\\
& = K^{n+1} + \frac{1}{n!} \sum_{k=1}^{n!} z_k^{n+1},
\end{align*}
and therefore we choose $K$ so that $K^{n+1} = A_{n+1}^{n+1} - \frac{1}{n!} \sum_{k=1}^{n!} z_k^{n+1}$. Using (\ref{estn}) we deduce
$|K|^{n+1} \leq A^{n+1} + D_n^{n+1} A^{n+1}$ and this gives $|K| \leq A (1 + D_n^{n+1})^{1/n+1}$. Finally  we get
$$|w_{k,l}| = |z_k + \xi^l K| \leq |z_k| + |K| \leq A + A (1 + D_n^{n+1})^{1/n+1} = D_{n+1} A$$
which completes the proof. $\Box$

By the above theorem we have $K_n \leq K_{n, n!} \leq D_n$. Let us prove that
\begin{equation}\label{lim}
\lim_{n\to\infty} \frac{D_n}{n} = 1.
\end{equation}
First, we have $D_{n+1} \geq 1 + D_n$ which gives $D_n \geq n$ for all $n \in \mathbb N$. Next
$$D_n = 1 + D_{n-1} (1 + D_{n-1}^{-n})^{1/n} \leq 1 + D_{n-1} \left( 1 + \frac{1}{nD_{n-1}^n} \right),$$
which implies
$$D_n - D_{n-1} \leq 1 + \frac{1}{nD_{n-1}^{n-1}}.$$
Since $\displaystyle\lim_{n \to \infty} D_n = +\infty$ we get, using Stolz Theorem on sequences
$$\limsup_{n \to \infty} \frac{D_n}{n} \leq \limsup_{n \to \infty} (D_n - D_{n-1}) \leq \limsup_{n\to\infty}
\left( 1 + \frac{1}{nD_{n-1}^{n-1}}\right) = 1,$$
which, combined with $D_n \geq n$, gives (\ref{lim}).

{\it Acknowledgment:} The authors are thankful to Du\v san Vojinovi\'c for programming which helped us to form conjectures
in the last part of this work.


\begin{thebibliography}{99}

\bibitem{Pr}   V, Prasolov, {\it Polynomials}, Springer-Verlag, 2004.

\bibitem{Ru} Z, Rubinstein, {\it Some results in the location of zeroes of polynomials}, Pacific J. of Mathematics
Vol. 15, No. 4, 1965.

\bibitem{Tc}   L. Tchakaloff, {\it Sur la distribution des zeros d'une classe des polynomes algebriques},
C. R. Acad. Bulgare Sci. 13 (1960), 249-251.


\end{thebibliography}
\end{document}